\documentclass{amsart}

 \usepackage[dvips]{graphicx}
 \usepackage{amssymb}
 \usepackage{mdwlist}

 \newtheorem{theorem}{Theorem}[section]
 \newtheorem{lemma}[theorem]{Lemma}
 \newtheorem{claim}[theorem]{Claim}
 \newtheorem{corollary}[theorem]{Corollary}

\theoremstyle{definition}
 \newtheorem{definition}[theorem]{Definition}
 \newtheorem{example}[theorem]{Example}

\theoremstyle{remark}
 \newtheorem{remark}[theorem]{Remark}

\begin{document}

\title[Lens spaces obtainable by surgery on knots] 
{Lens spaces obtainable by surgery\\ on doubly primitive knots}

\author{Kazuhiro Ichihara} 
\thanks{The first author is supported in part by Grant-in-Aid for Young Scientists (No. 18740038), The Ministry of Education, Culture, Sports, Science and Technology, Japan. }
\address{School of Mathematics Education, 
Nara University of Education, 
Takabatake-cho, Nara 630-8528, Japan.} 
\email{ichihara@nara-edu.ac.jp}

\author{Toshio Saito} 
\thanks{The second author is partially supported by 
JSPS Research Fellowships for Young Scientists.}
\address{Graduate School of Humanities and Sciences,
Nara Women's University,
Kitauoyanishi-machi, Nara 630--8506, Japan.} 
\email{tsaito@cc.nara-wu.ac.jp}



\date{\today}

\begin{abstract} 
In this paper, we consider which lens spaces are obtainable by 
Dehn surgery described by Berge on doubly primitive knots. 
It is given an algorithm to decide 
whether a given lens space is obtainable by such surgery. 
Also included is a complete characterization of 
such surgery yielding lens spaces with Klein bottles. 
\end{abstract} 
\maketitle

\section{Introduction}

Given a knot $K$ in a $3$-manifold, 
the following operation is called \textit{Dehn surgery} on $K$: 
Remove an open regular neighborhood of $K$, and glue a solid torus back. 
Dehn surgeries on the trivial knot in the 3-sphere $S^3$ give the well-known class of 3-manifolds, so-called \textit{lens spaces}. 
On the other hand, it is known that only restricted Dehn surgeries on non-trivial knots yield lens spaces. 
Such a non-trivial example was found by Fintushel and Stern in \cite{FS}, and 
then a lot of researches have been done about them. 

Based on the pioneering works by Berge \cite{Berge}, 
Gordon conjectured in \cite[Problem 1.78]{K} that 
the knots admitting Dehn surgeries yielding lens spaces are all \textit{doubly primitive}. 
In this paper, we concentrate our attention to such knots, and 
consider the problem: 
Which lens spaces are obtainable by Dehn surgery on such knots? 
Please see the next section for the definitions in detail.

We first give an algorithm to decide whether a given lens space is obtainable 
by Dehn surgeries described by Berge on doubly-primitive knots. 
Our algorithm depend on the works on the triviality of three bridge knots presented in \cite{HO}, 
and so, it is quite effective. 

Next we consider a particular class of lens spaces; the ones containing Klein bottles. 
In the study of Dehn surgeries giving 3-manifolds with Klein bottles, 
Teragaito asked the following question: 
Which non-trivial knots admit Dehn surgery yielding lens spaces with Klein bottles? 
In particular, can hyperbolic knots admit such surgeries? 
Our result gives a partial answer to the question:

\begin{theorem}\label{main theorem}
Dehn surgery on a doubly primitive knot described by Berge 
yields a lens space containing a Klein bottle 
if and only if the knot is either of the $(\pm 5,3)$-, or, 
$(\pm 7,3)$-torus knots. The lens spaces so-obtained are of type $(16,7)$ or $(20,9)$. 
In particular, no hyperbolic doubly primitive knots admit such Dehn surgeries. 
\end{theorem} 

We remark that 
a characterization of lens spaces containing Klein bottles 
was already established in \cite{BW}. 

Also remark that it has already known completely 
which lens spaces are obtainable by surgery on torus knots in \cite{Moser}. 
Furthermore, by \cite{BL,Wang,Wu}, 
it was shown that 
satellite knots admitting Dehn surgery yielding lens spaces are 
all certain cabled knots, which are shown to be doubly primitive. 
Thus our theorem implies the following corollary immediately. 

\begin{corollary}
Dehn surgery on a non-trivial non-hyperbolic knot 
yields a lens space containing a Klein bottle 
if and only if the knot is either of the $(\pm 5,3)$-, or, 
$(\pm 7,3)$-torus knots. The lens spaces so-obtained are of type $(16,7)$ or $(20,9)$. 
In particular, no satellite knots admit Dehn surgery 
yielding a lens space containing a Klein bottle. 
\end{corollary}

About the question above, recently, Tange showed that 
only lens spaces of type $(16,7)$ or $(20,9)$ are obtainable by Dehn surgeries on 
non-trivial knots. However nothing could be said there about the types of the knots 
admitting such surgeries.

\section{Preliminaries} 
In this section, we will set up our terminologies. 
In the following, $E(B;A)$ denotes 
the \textit{exterior} of a subset $B$ in a topological space $A$, 
i.e., $E(B;A)=\mathrm{cl}(A\backslash \eta(B;A))$, 
where $\eta(B;A)$ means an open regular neighborhood of $B$ in $A$. 

\begin{figure}[htb]\begin{center}
  {\unitlength=1cm
  \begin{picture}(12,4.2)
  \put(3.05,0.4){\includegraphics[keepaspectratio]
  {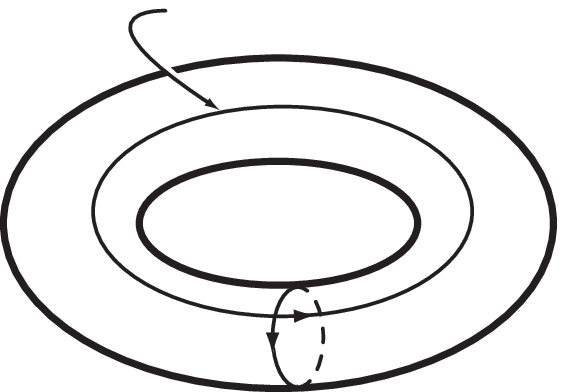}}
   \put(4.9,4.2){$\ell$}
   \put(5.9,0){$m$}
  \end{picture}}
  \caption{}
  \label{Fig1}
\end{center}\end{figure}

\subsection{Lens space}
Let $V_1$ be a standard solid torus in $S^3$, $m$ a meridian of $V_1$ and $\ell$ a 
longitude of $V_1$ such that $\ell$ bounds a disk in $\mathrm{cl}(S^3\setminus V_1)$. 
We fix an orientation of $m$ and $\ell$ as illustrated in Figure \ref{Fig1}. 
Let $p$ and $q$ be coprime integers and $\bar{m}$ a meridian of $V_2$. 
Then, by attaching another solid torus $V_2$ to $V_1$ so that $\bar{m}$ is isotopic to a representative 
of the homology class $p[\ell]+q[m]$, 
we obtain a \textit{lens space}, denoted by $L(p,q)$.

\subsection{Dehn surgery}
Let $K$ be a knot in a connected compact orientable $3$-manifold $N$. 
We fix an oriented meridian-longitude system $(m,\ell)$ for $K$ like in Figure \ref{Fig1}. 
When $K\subset S^3$, we always take the preferred longitude for $K$ as $\ell$. 
Recall that a \textit{Dehn surgery} on a knot $K$ is an operation to attach a solid torus 
$\bar{V}$ to $E(K;N)$ by a boundary-homeomorphism $\varphi:\partial \bar{V}\to \partial E(K;N)$. 
If $\varphi(\bar{m})$ is isotopic to a representative of the homology class $p[m]+q[\ell]$ for a meridian $\bar{m}$ of $\bar{V}$, 
then the surgery is called \textit{$p/q$-surgery}. 
By an \textit{integral surgery}, we mean an $r$-Dehn surgery with $r$ an integer. 
Set $N_{\varphi}=E(K;N)\cup_{\varphi} \bar{V}$ and let $K^{\ast}\subset N_{\varphi}$ be a core loop of $\bar{V}$. 
We call $K^{\ast}$ the \textit{dual knot} of $K$ in $N_{\varphi}$. 
We remark that $E(K;N)$ is homeomorphic to $E(K^{\ast};N_{\varphi})$ and that 
if a Dehn surgery on $K$ in $N$ yields a $3$-manifold $N_{\varphi}$, 
then $K^{\ast}$ admits a Dehn surgery yielding $N$.

\subsection{Doubly primitive knot and dual knot}
Let $H$ be a genus two handlebody standardly embedded in $S^3$, 
i.e., $E(H;S^3)$ is also a genus two handlebody. 
A simple closed curve on the boundary $\partial H$ is in a \textit{doubly primitive position} 
if it represents a free generator both of $\pi_1(H)$ and of $\pi_1(E(H;S^3))$. 
A knot in $S^3$ is called a \textit{doubly primitive knot} 
if it is isotopic to a simple closed curve in a doubly primitive position. 
Let $K$ be a doubly primitive knot. 
When $K$ is isotoped into a doubly primitive position, 
$\partial H\cap  \partial \eta(K;S^3)$ consists of two essential simple closed curves  
which are mutually isotopic in $\partial \eta(K;S^3)$. 
The isotopy class is called a \textit{surface slope} of $K$. 
We remark that a Dehn surgery along a surface slope of $K$ is always an integral surgery. 
It is then proved by Berge \cite{Berge} that any Dehn surgery along a surface slope of $K$ 
yields a \textit{lens space} $L(p,q)$, which is also obtained by $(-q/p)$-surgery on a trivial knot, 
for some integers $p$ and $q$. 

Moreover, he showed the dual knot of $K$ in $L(p,q)$ is isotopic to a knot defined as follows. 

\begin{figure}[htb]\begin{center}
  {\unitlength=1cm
  \begin{picture}(12,5.52)
  \put(2.41,0){\includegraphics[keepaspectratio]
  {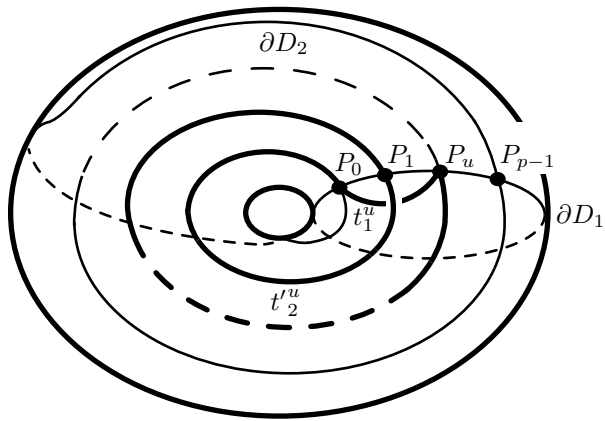}}
   \put(5.75,4.9){$\partial D_2$}
   \put(6.75,3.3){$P_0$}
   \put(7.45,3.4){$P_1$}
   \put(8.22,3.4){$P_u$}
   \put(8.95,3.4){$P_{p-1}$}
   \put(7,2.58){$t^u_1$}
   \put(9.7,2.6){$\partial D_1$}
   \put(5.9,1.45){${t'}^u_2$}
  \end{picture}}
  \caption{Here, ${t'}^u_2$ is a projection of $t^u_2$ on $\partial V_1$. }
  \label{Fig2}
\end{center}\end{figure}

\begin{definition}\label{def}
Let $V_1$ be a standard solid torus in $S^3$, $m$ a meridian of $V_1$ and 
$\ell$ a longitude of $V_1$ such that $\ell$ bounds a disk in $\mathrm{cl}(S^3\setminus V_1)$. 
We fix an orientation of $m$ and $\ell$ as illustrated in Figure \ref{Fig1}. 
By attaching a solid torus $V_2$ to $V_1$ so that $\bar{m}$ is isotopic to a representative of $p[\ell]+q[m]$, 
we obtain a lens space $L(p,q)$, where $p$ and $q$ are coprime integers and $\bar{m}$ is a meridian of $V_2$. 
The intersection points of $m$ and $\bar{m}$ are labeled by $P_0,\ldots ,P_{p-1}$ 
successively along the positive direction of $m$. 
Let $t^u_i$ $(i=1,2)$ be simple arcs in $D_i$ joining $P_0$ to $P_u$ $(u=1,2,\dots,p-1)$. 
Then the notation $K(L(p,q);u)$ denotes the knot $t^u_1\cup t^u_2$ in $L(p,q)$. 
See Figure \ref{Fig2}. 
\end{definition}

If $K(L(p,q);u)$ is the dual knot of some doubly primitive knot in $S^3$, 
then it admits a Dehn surgery yielding $S^3$. 
This Dehn surgery is \textit{longitudinal}, that is, 
it is an integral surgery for a meridian-longitude system. 
We remark that the converse does not hold in general; 
it is not necessary for any knot represented by $K(L(p,q);u)$ 
to admit a longitudinal surgery yielding $S^3$.

\begin{remark}\label{remark}
It is known that two lens spaces $L(p,q)$ and $L(p',q')$ are 
(possibly orientation reversing) homeomorphic if and only if 
$|p|=|p'|$, and $q\equiv \pm q'$ $(\bmod\ p)$ or $qq'\equiv \pm 1$ $(\bmod\ p)$.  
Also, we easily see that $K(L(p,q);u)$ is isotopic to 
$K(L(p,q);p-u)$. Hence for $K(L(p,q);u)$, we assume that $0<q<p/2$ and 
$1\le u\le p/2$ in the remainder of the paper. 
\end{remark}

\section{Algorithm to detect obtainable lens spaces} 
In this section, we will describe an algorithm to decide whether a given 
lens space is obtainable by Dehn surgeries on doubly-primitive knots along surface slopes. 

As an instructive example, let us check that the knot $K(L(5,1);2)$ in $L(5,1)$ can admit 
Dehn surgery creating $S^3$. 
This implies that the lens space $L(5,1)$ is obtainable by Dehn surgery on some doubly primitive 
knot in $S^3$. 

\begin{figure}[htb]\begin{center}
  {\unitlength=1cm
  \begin{picture}(12,7.9)
  \put(0.17,0.3){\includegraphics[keepaspectratio]
  {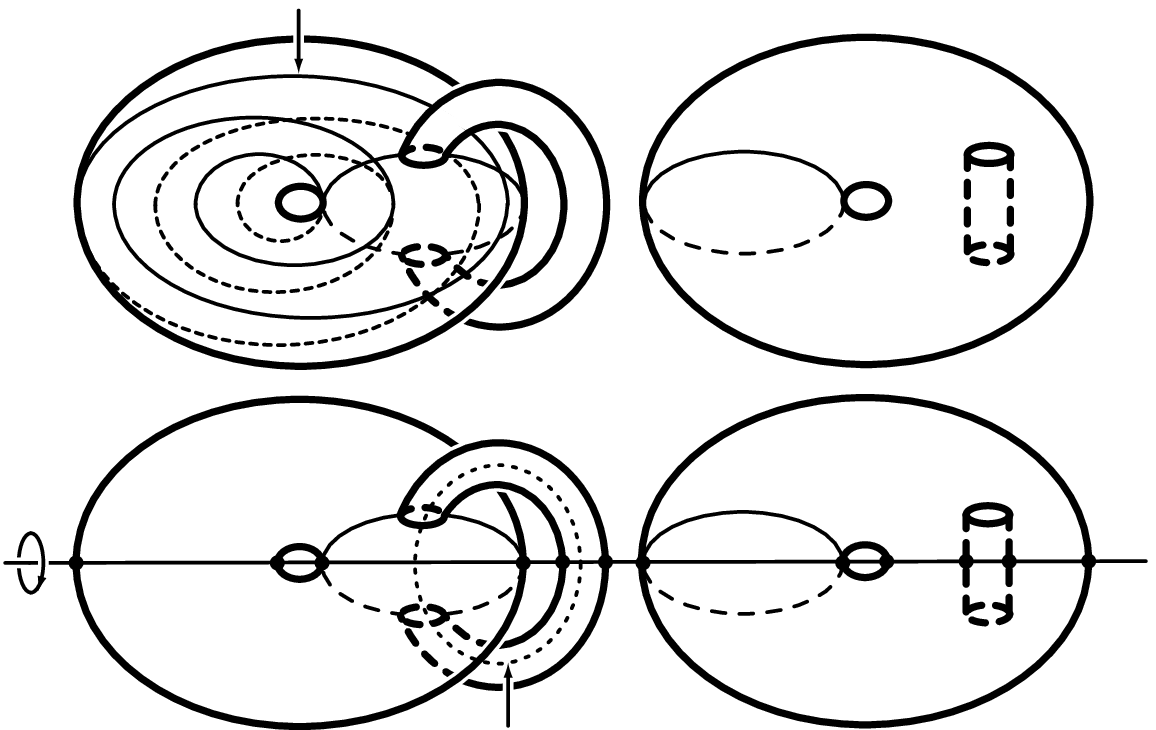}}
   \put(3,7.75){$\partial D_2$}
   \put(7.5,6.3){$D_2$}
   \put(5.2,0){$K(L(5,1);2)$}
  \end{picture}}
  \caption{}
  \label{Fig3}
\end{center}\end{figure}

\begin{figure}[htb]\begin{center}
  {\unitlength=1cm
  \begin{picture}(12,3)
  \put(0.78,0){\includegraphics[keepaspectratio]
  {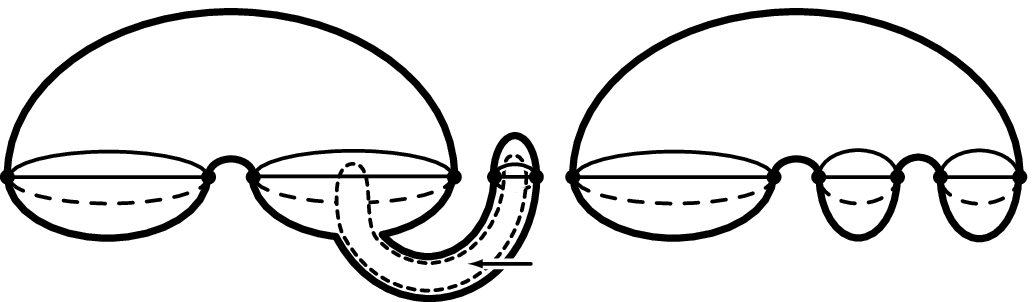}}
   \put(6.25,0.3){$B$}
  \end{picture}}
  \caption{}
  \label{Fig4}
\end{center}\end{figure}

\begin{figure}[htb]\begin{center}
  {\unitlength=1cm
  \begin{picture}(12,3.72)
  \put(3.25,0){\includegraphics[keepaspectratio]
  {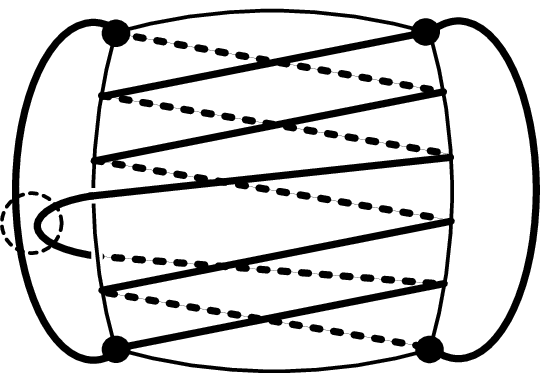}}
   \put(3.55,1.9){$B$}
  \end{picture}}
  \caption{}
  \label{Fig5}
\end{center}\end{figure}

\begin{figure}[htb]\begin{center}
  {\unitlength=1cm
  \begin{picture}(12,3.7)
  \put(0.29,0){\includegraphics[keepaspectratio]
  {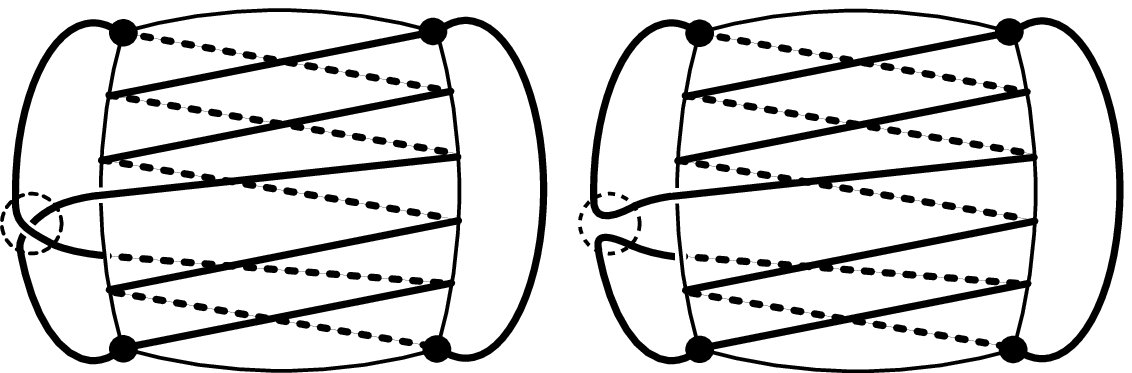}}
   \put(0.6,1.9){$B$}
   \put(6.45,1.9){$B$}
  \end{picture}}
  \caption{}
  \label{Fig6}
\end{center}\end{figure}

\begin{example}
We consider a Heegaard splitting of genus two of $L(5,1)$ illustrated as in Figure \ref{Fig3}, 
which is obtained from the standard Heegaard splitting of genus one by stabilization. 
Note that the knot $K(L(5,1);2)$ is isotopic to the dotted knot in Figure \ref{Fig3} below. 
Take the quotient of $L(5,1)$ by involution as illustrated in Figure \ref{Fig3}. 
It follows from Figure \ref{Fig4} that the quotient space is $S^3$. 
Let $B$ be the $3$-ball in Figure \ref{Fig4}, 
which appears as the quotient of a equivariant regular neighborhood of $K(L(5,1);2)$ in $L(5,1)$. 
Also the quotient of the axis of the involution gives a knot in $S^3$ as shown in Figure \ref{Fig5}.

If $K(L(5,1);2)$ admits a Dehn surgery yielding $S^3$, 
then, from the knot shown in Figure \ref{Fig5}, 
the corresponding untangle surgery at the $3$-ball $B$ must give the trivial knot. 
This follows from the so-called \textit{Montesinos trick} originally developed in \cite{Montesinos}. 

In fact, it suffice to check the only two links shown in Figure \ref{Fig6}, 
due the result of the second author given in \cite[Theorem 2.5]{Saito}. 
We can see that the knot depicted in the left side of Figure \ref{Fig6} is actually trivial. 
Therefore $K(L(5,1);2)$ can admit Dehn surgery yielding $S^3$. 
\end{example}

In general, by the following procedure, we can determine whether the given $L(p,q)$ is obtainable by 
a Dehn surgery on a non-trivial doubly primitive knot along a surface slope or not. 

\begin{figure}[htb]\begin{center}
  {\unitlength=1cm
  \begin{picture}(12,4.18)
  \put(2.92,0){\includegraphics[keepaspectratio]
  {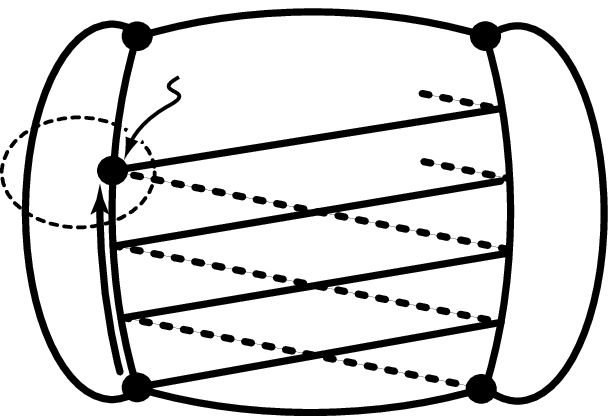}}
   \put(3.55,3.15){$B$}
   \put(4.77,3.45){the $u$-th vertex $V$}
  \end{picture}}
  \caption{}
  \label{Fig7}
\end{center}\end{figure}

\begin{figure}[htb]\begin{center}
  {\unitlength=1cm
  \begin{picture}(12,3.7)
  \put(0.29,0){\includegraphics[keepaspectratio]
  {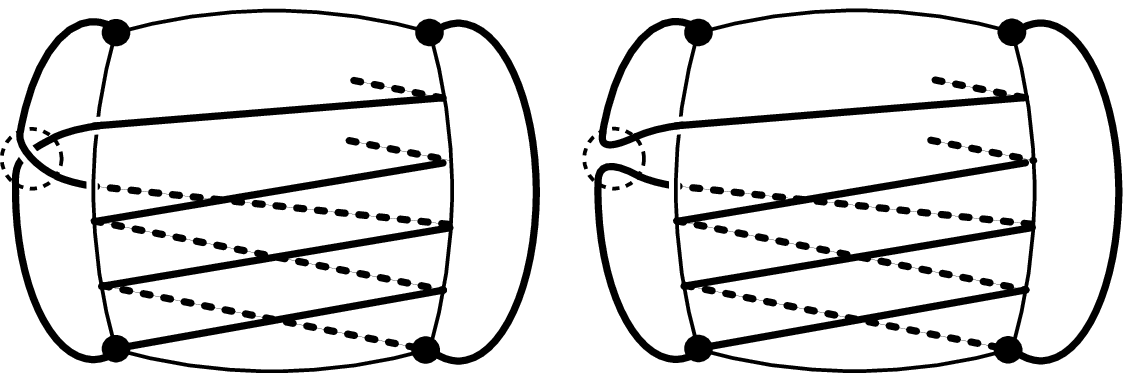}}
   \put(0.65,2.6){$B$}
   \put(6.55,2.6){$B$}
  \end{picture}}
  \caption{}
  \label{Fig8}
\end{center}\end{figure}

\begin{enumerate}
\item 
Consider the two-bridge link of type $(p,q)$ represented by the Schubert form $\mathrm{b}(p,q)$. 
This knot has the diagram illustrated as in Figure \ref{Fig7}. 
See \cite{BZ} for example. 

\suspend{enumerate}

For any integer $u$ with $1\le u\le p/2$, we do the following steps $(2)$-$(4)$ repeatedly. 
\resume{enumerate}

\item 
Put a vertex $V$ on the diagram as illustrated in Figure \ref{Fig7}, that is, 
put $V$ on the $u$-th ``wedge'' from the left-bottom side of the pillowcase. 

\item 
In the neighborhood of $V$, depicted as the encircled region in Figure \ref{Fig7}, 
make a crossing (the left) or do smoothing (the right) as in Figure \ref{Fig8}. 

\item 
Exactly one of the two diagrams so obtained gives a $3$-bridge knot (not a link). 
By the algorithm given in \cite{HO}, we check the knot is trivial or not. 
The knot is trivial if and only if $K(L(p,q);u)$ admits a Dehn surgery yielding $S^3$, that is, 
$L(p,q)$ is obtainable by 
a Dehn surgery on a non-trivial doubly primitive knot along a surface slope. 
\end{enumerate}

\section{Proof of Theorem \ref{main theorem}} 

In this section, we consider the lens spaces containing Klein bottles, 
and give a proof of Theorem \ref{main theorem}. 
In this case, by virtue of the result \cite[Theorem 2.5]{Saito} of the second author, 
we can determine the possible position of the vertex $V$ 
in the algorithm in the previous section. 
To state \cite[Theorem 2.5]{Saito}, we prepare the following notations. 

\begin{definition} 
Let $p$ and $q$ be a pair of positive coprime integers. 
Let $\{s_j\}_{1\le j\le p}$ be the finite sequence such that $0\le s_j< p$ and $s_j\equiv q\cdot j$ ($\bmod\ p$). 
For an integer $k$ with $0<k<p$, 
$\Psi_{p,q} (k)$ denotes the smallest integer $j$ with $s_j=k$, and 
$\Phi_{p,q} (k)$ the number of elements of the set 
$\{s_j|\ 1\le j< \Psi_{p,q} (k),\ s_j< k\}$. 
\end{definition}

Then the following is our key lemma. 

\begin{lemma}[{\cite[Theorem 2.5]{Saito}}]\label{TS1}\label{formula}
Let $p$ and $q$ be coprime integers with $0<q<p$ and $u$ an integer with $1\le u\le p-1$. 
If $K(L(p,q);u)$ admits longitudinal surgery yielding $S^3$, then we have 

\medskip\hspace{2.5cm}
$p\cdot \Phi_{p,q} (u)-u\cdot \Psi_{p,q} (u)=\pm 1$ or $\pm 1-p$. 

\medskip
\noindent
In particular, $p$ and $u$ are coprime. 
\end{lemma}

By using this, we prove our theorem as follows.

\begin{proof}[Proof of Theorem \ref{main theorem}]
Let $K$ be a doubly primitive knot with a Dehn surgery along a surface slope yielding a lens space $M$ containing a Klein bottle. 
Then $M$ have to be of type $(4n,2n-1)$, where $n$ is a positive integer, 
up to (possibly orientation-reversing) homeomorphism \cite{BW}. 
Hence the dual knot $K^{\ast}$ of $K$ in $M$ 
is represented by $K(L(4n,2n-1);u)$ for some positive integer $u$. 
Note that this $K^{\ast}$ have to admit a longitudinal surgery yielding $S^3$.

\begin{claim}\label{main prop}
If $K(L(4n,2n-1);u)$ with $(n\ge 2)$ admits 
longitudinal surgery yielding $S^3$, then $\Phi_{4n,2n-1} (u)=0$ and 
one of the following holds. 

\begin{center}
$(n,u)=(4,3),(4,5),(5,3),(5,7)$
\end{center}
\end{claim}

\begin{proof}
Since $K(L(4n,2n-1);u)$ admits a longitudinal surgery yielding $S^3$, 
it follows from Lemma \ref{formula} that $4n$ and $u$ are coprime. 
Hence $u$ is an odd integer. 
By Remark \ref{remark}, we can assume that $u<2n$. 
Then the sequence $\{u_j\}_{1\le j\le 4n}$ for $(4n,2n-1)$ is : 

\begin{center}
$u_{j}\equiv\left\{
       \begin{array}{ll}
       2n-j\ (\bmod\ p) & \mathrm{if}\ j\ \mathrm{is\ odd}\\
       4n-j\ (\bmod\ p) & \mathrm{if}\ j\ \mathrm{is\ even}  
       \end{array}
\right.$. 
\end{center}

In particular, a sub-sequence $\{u_j\}_{1\le j\le 2n-1}$ of 
$\{u_j\}_{1\le j\le 4n}$ satisfies the following. 

\begin{enumerate}
\item
$u_j$ is odd if $j$ is odd, and $u_j$ is even if $j$ is even. 
\item
Each of sub-sequences $\{u_{2k-1}\}_{1\le k\le n}$ and 
$\{u_{2k}\}_{1\le k\le n}$ is monotonically decreasing. 
\item
$\max\{u_{2k-1}\ |\ \ 1\le k\le n\}=u_1=2n-1$ and $\min\{u_{2k}\ |\ \ 1\le k\le n\}=u_{2n-2}=2n-2$. 
Hence we have $\max\{u_{2k-1}\ |\ \ 1\le k\le n\}<\min\{u_{2k}\ |\ \ 1\le k\le n\}$. 
\end{enumerate}

Let $m$ be the integer satisfying $u_m=u$, that is, $m=\Psi_{4n,2n-1} (u)$. 
Then $u=2n-m$. Since $u$ is an odd integer less that $2n$, we see that $1\le m\le 2n-1$. 
This implies that $u_j>u$ for any integer $j$ satisfying $1\le j\le m-1$ and hence $\Phi_{4n,2n-1} (u)=0$. 
Therefore we have the first conclusion of the claim. 

By Lemma \ref{formula}, we also have : 

\begin{center}
$u\cdot m=\pm 1$ or $4n\pm 1$.
\end{center}

\medskip\noindent
\textit{Case 1.}\ \ $u\cdot m=\pm 1$.

\medskip
Since $u$ and $m$ are positive integers, we have $u=m=1$. This implies that 
$u_1=1$ and hence $n=1$. This contradicts that $n\ge 2$. 

\medskip\noindent
\textit{Case 2.}\ \ $u\cdot m=4n\pm 1$.

\medskip
In this case, we first note that $m\ne 2$. Hence we have the following.  

\begin{eqnarray*}
u\cdot m      &=& 4n\pm 1\\
(2n-m)\cdot m &=& 4n\pm 1\\
2n            &=& m+2+\displaystyle\frac{\>4\pm 1\>}{m-2}
\end{eqnarray*}

Since $m$ and $n$ are positive integers, we have the desired conclusion. 
\end{proof}

By this claim, the possible type of the obtained lens space $M$ is $(16,7)$ or $(20,9)$. 

In fact, by applying our algorithm, we can directly check that 
$K(L(4n,2n-1);u)$ actually admits longitudinal surgery yielding $S^3$
for $(n,u)=(4,3),(4,5),(5,3),(5,7)$.

Let us consider these cases in detail. 

\medskip\noindent
\textit{Case 1.}\ \ $K^{\ast}=K(L(16,7);3)$ or $K(L(16,7);5)$. 

\medskip
Since $\Phi_{16,7} (3)=\Phi_{16,7} (3)=0$, 
we can see that $K^{\ast}$ is a torus knot (\textrm{cf}. \cite[Proposition 5.2]{Saito2}). 
Since $E(K^{\ast};M)\cong E(K;S^3)$, $E(K;S^3)$ admits 
Seifert fibration and hence $K$ is a torus knot. 
It then follows from Van-Kampen's theorem (\textrm{cf}. \cite[Section 5]{Saito3}) that 

\begin{eqnarray*}
\pi_1 (E(K;S^3)) 
&\cong & \pi_1 (E(K^{\ast};M))\\
&\cong & \langle x,y\ |\ x^5=y^3 \rangle.
\end{eqnarray*}

This implies that $K$ is the $(\pm 5,3)$-torus knot, and 
the obtained lens space $M$ is homeomorphic to $L(16,7)$. 

\medskip\noindent
\textit{Case 2.}\ \ $K^{\ast}=K(L(20,9);3)$ or $K(L(20,9);7)$. 

\medskip
In the same way as above, we see that $K$ is a torus knot. 
It also follows from Van-Kampen's theorem that 

\begin{eqnarray*}
\pi_1 (E(K;S^3)) 
&\cong & \pi_1 (E(K^{\ast};M))\\
&\cong & \langle x,y\ |\ x^7=y^3 \rangle.
\end{eqnarray*}

This implies that $K$ is the $(\pm 7,3)$-torus knot, and 
the obtained lens space $M$ is homeomorphic to $L(20,9)$. 
\end{proof}



\begin{thebibliography}{99} 

\bibitem{Berge}J. Berge, 
Some knots with surgeries yielding lens spaces, 
unpublished manuscript. 

\bibitem{BL}S. Bleiler and R. Litherland, 
Lens spaces and Dehn surgery, 
Proc. Amer. Math. Soc. \textbf{107} (1989), 1127--1131. 

\bibitem{BW}
G. E. Bredon and J. W. Wood, 
Non-orientable surfaces in orientable $3$-manifolds, 
Invent. Math. \textbf{7} (1969), 83--110. 

\bibitem{BZ}G. Burde and H. Zieschang, 
Knots, 
de Gruyter, Berlin and New York (1986).

\bibitem{FS}
R. Fintushel and R. J. Stern, 
Constructing lens spaces by surgery on knots, 
Math. Z. \textbf{175} (1980), 33--51. 

\bibitem{HO}T. Homma and M. Ochiai, 
On relations of Heegaard diagrams and knots, 
Math. Sem. Notes Kobe Univ. \textbf{6} (1978), 383--393. 

\bibitem{K} R. Kirby, 
Problems in low-dimensional topology, 
Geometric topology, AMS/IP Stud. Adv. Math., 2.2, (Athens, GA, 1993), 
(Amer. Math. Soc., Providence, RI, 1997), 35--473. 

\bibitem{Montesinos}Jos\'{e} M. Montesinos, 
Surgery on links and double branched covers of $S\sp{3}$. 
Knots, groups, and $3$-manifolds (Papers dedicated to the memory of R. H. Fox), 
pp. 227--259. Ann. of Math. Studies, No. 84, Princeton Univ. Press, 
Princeton, N.J., 1975. 

\bibitem{Moser}L. Moser, 
Elementary surgery along a torus knot, 
Pacific J. Math. \textbf{38} (1971), 734--745. 

\bibitem{Saito}T. Saito, 
Dehn surgery and (1,1)-knots in lens spaces, 
Topology Appl. \textbf{154} (2007), no. 7, 1502--1515.

\bibitem{Saito2}T. Saito, 
The dual knots of doubly primitive knots, 
Osaka J. Math. (to appear). 

\bibitem{Saito3}T. Saito, 
Knots in lens spaces with the 3-sphere surgery, preprint. 

\bibitem{Wang}S. Wang, 
Cyclic surgery on knots, 
Proc. Amer. Math. Soc. \textbf{107} (1989), 1091--1094. 

\bibitem{Wu}Y. Q. Wu, 
Cyclic surgery and satellite knots, 
Topology Appl. \textbf{36} (1990), 205--208. 
\end{thebibliography}
\end{document}